\documentclass[12pt]{amsart}
\usepackage{times}

\newtheorem{thm}{Theorem}
\newtheorem{lemma}[thm]{Lemma}

\newcommand{\D}{{\mathbb D}}
\newcommand{\C}{{\mathbb C}}
\newcommand{\bloch}{{\mathcal B}}
\newcommand{\cn}{\C^n}
\newcommand{\bn}{{\mathbb B}_n}

\newcommand{\ind}{\int_\D}
\newcommand{\daa}{dA_\alpha}

\begin{document}

\title[Bergman Spaces]
{Lipschitz Type Characterizations\\
 for Bergman Spaces}

\author{Hasi Wulan}
\address{Department of Mathematics\\
         Shantou University\\
         Guangdong, China}
\email{wulan@stu.edu.cn}

\author{Kehe Zhu}
\address{Department of Mathematics\\
              SUNY\\
              Albany, NY 12222, USA\\
              and Department of Mathematics\\
              Shantou University\\
              Guangdong, China}
\email{kzhu@math.albany.edu}

\thanks{Research partially supported by the National Science Foundations 
of the US and China}

\begin{abstract}
We obtain new characterizations for Bergman spaces with standard
weights in terms of Lipschitz type conditions in the Euclidean, 
hyperbolic, and pseudo-hyperbolic metrics. As a consequence, we
prove optimal embedding theorems when an analytic function
on the unit disk is symmetrically lifted to the bidisk.
\end{abstract}

\maketitle

\section{Introduction}

Let $\D$ be the open unit disk in the complex plane $\C$. For any
$\alpha>-1$ we consider the weighted area measure
$$dA_\alpha(z)=(\alpha+1)(1-|z|^2)^\alpha\,dA(z),$$
where $dA$ is the normalized area measure on $\D$. It is easy to see 
that each $dA_\alpha$ is a probability measure on $\D$.

For $p>0$ and $\alpha>-1$ we denote by $A^p_\alpha$ the space of analytic
functions $f$ in $\D$ such that
$$\ind|f(z)|^p\,\daa(z)<\infty.$$
These are called weighted Bergman spaces with standard weights. See 
\cite{DS} and \cite{HKZ} for the modern theory of Bergman spaces.

Three different metrics on the unit disk will be used in the paper.
First, the usual Euclidean metric is of course written as $|z-w|$.
Second, the pseudo-hyperbolic metric on $\D$ is given by
$$\rho(z,w)=\left|\frac{z-w}{1-\overline zw}\right|.$$
And finally, the hyperbolic metric on $\D$ is given by
$$\beta(z,w)=\frac12\log\frac{1+\rho(z,w)}{1-\rho(z,w)}.$$
We mention that the hyperbolic metric is also called the Bergman metric,
and sometimes the Poincare metric as well.

The main result of the paper is the following.

\begin{thm}
Suppose $p>0$, $\alpha>-1$, and $f$ is analytic in $\D$. Then the
following conditions are equivalent.
\begin{itemize}
\item[(a)] $f$ belongs to $A^p_\alpha$.
\item[(b)] There exists a continuous function $g$ in $L^p(\D,\daa)$ such that
$$|f(z)-f(w)|\le\rho(z,w)(g(z)+g(w))$$
for all $z$ and $w$ in $\D$.
\item[(c)] There exists a continuous function $g$ in $L^p(\D,\daa)$ such that
$$|f(z)-f(w)|\le\beta(z,w)(g(z)+g(w))$$
for all $z$ and $w$ in $\D$.
\item[(d)] There exists a continuous function $g$ in
$L^p(\D,dA_{p+\alpha})$ such that
$$|f(z)-f(w)|\le|z-w|(g(z)+g(w))$$
for all $z$ and $w$ in $\D$.
\end{itemize}
\label{1}
\end{thm}

Note that the same measure $dA_\alpha$ appears in conditions (a), (b), and
(c), but condition (d) involves a different measure $dA_{p+\alpha}$.

Similar characterizations for Hardy-Sobolev spaces have appeared in the
literature before. See \cite{H1}\cite{H2} for example. The present
paper was motivated by \cite{S}. As another motivation for our result, we mention
that the classical Bloch space $\bloch$, consisting of analytic functions $f$ in
$\D$ such that
$$\sup\{(1-|z|^2)|f'(z)|:z\in\D\}<\infty,$$
also admits a Lipschitz type characterization. More specifically, an analytic function
$f$ in $\D$ belongs to $\bloch$ if and only if there exists a positive constant $C$
such that
$$|f(z)-f(w)|\le C\beta(z,w)$$
for all $z$ and $w$ in $\D$; see \cite{HKZ} or \cite{Z} for example. It is then clear
that an analytic function $f$ in $\D$ belongs to $\bloch$ if and only if there exists
a continuous function $g$ in $L^\infty(\D)$ such that
$$|f(z)-f(w)|\le\beta(z,w)(g(z)+g(w))$$
for all $z$ and $w$ in $\D$.

Since the Bergman metric is based on the reproducing kernel of the Bergman 
space, it is no surprise that our Lipschitz type characterizations for weighted
Bergman spaces appear more natural when the Bergman metric (and its bounded
counterpart, the pseudo-hyperbolic metric) is used. Althought a characterization
in terms of the Euclidean metric (condition (d)) is possible, it gives one the 
impression of something artificial.

\section{Preliminaries}

We collect some preliminary results in this section that involve the
hyperbolic and pseudo-hyperbolic metrics.

For any $0<r<1$ and $z\in\D$ we let
$$D(z,r)=\{w\in\D:\rho(w,z)<r\}$$
denote the pseudo-hyperbolic disk centered at $z$ with radius $r$. It is
well known that $D(z,r)$ is actually a Euclidean disk with Euclidean center
and Euclidean radius given by
$$\frac{1-r^2}{1-r^2|z|^2}\,z,\qquad\frac{1-|z|^2}{1-r^2|z|^2}\,r,$$
respectively; see \cite{G} for example. In particular, if $r$ is fixed, then
the area of $D(z,r)$ is comparable to $(1-|z|^2)^2$.

The pseudo-hyperbolic metric is bounded by $1$. But the hyperbolic metric
is unbounded. For any $R>0$ and $z\in\D$ we let
$$E(z,R)=\{w\in\D:\beta(w,z)<R\}$$
denote the hyperbolic disk centered at $z$ with radius $R$. If $0<r<1$ and
\begin{equation}
R=\frac12\log\frac{1+r}{1-r},
\label{eq1}
\end{equation}
then it is clear that $E(z,R)=D(z,r)$. Consequently, if $R$ is fixed, then
the area of $E(z,R)$ is comparable to $(1-|z|^2)^2$ as well. By the same
token, any estimate in terms of the pseudo-hyperbolic metric can be translated
to one in terms of the hyperbolic metric, and vise versa.

\begin{lemma}
For any fixed $r\in(0,1)$ there exists a positive constant $C$ such that
$$C^{-1}\le\frac{1-|z|^2}{|1-\overline zw|}\le C$$
whenever $\rho(z,w)\le r$. Consequently, there exists a positive constant
$C$ such that
$$C^{-1}\le\frac{1-|z|^2}{1-|w|^2}\le C$$
whenever $\rho(z,w)\le r$.
\label{2}
\end{lemma}

\begin{proof}
This is well known. See \cite{HKZ} or \cite{Z} for example.
\end{proof}

\begin{lemma}
For any fixed $r\in(0,1)$ and there exists a positive constant $C$
such that
$$|f(z)|^p\le\frac C{(1-|z|^2)^2}\int_{D(z,r)}|f(w)|^p\,dA(w)$$
for all $z\in\D$, all $p>0$, and all analytic $f$ in $\D$.
\label{3}
\end{lemma}

\begin{proof}
This is well known as well. See \cite{HKZ} or \cite{Z} for example.
\end{proof}

\begin{lemma}
For any fixed $z\in\D$ we have
$$\lim_{w\to z}\frac{\beta(z,w)}{|z-w|}=\lim_{w\to z}\frac{\rho(z,w)}{|z-w|}
=\frac1{1-|z|^2}.$$
\label{4}
\end{lemma}

\begin{proof}
This follows from elementary calculations.
\end{proof}

\begin{lemma}
For any $\alpha>-1$ and $p>0$ there exists a constant $C>0$ such that
$$\ind|f(z)|^p\,\daa(z)\le C\left[|f(0)|^p+\ind|g(z)|^p\,\daa(z)\right]$$
and
$$|f(0)|^p+\ind|g(z)|^p\,\daa(z)\le C\ind|f(z)|^p\,\daa(z)$$
for all analytic functions $f$ in $\D$, where
$$g(z)=(1-|z|^2)f'(z)$$
for $z\in\D$.
\label{5}
\end{lemma}

\begin{proof}
See \cite{HKZ} or \cite{Z} for example.
\end{proof}

\section{The Main Result}

We now prove the main result of the paper. The proof is constructive in
the sense that we will actually produce a formula for the function $g$
that appears in various conditions of Theorem~\ref{1}. We prove 
Theorem~\ref{1} as three separate results.

\begin{thm}
Suppose $p>0$, $\alpha>-1$, and $f$ is analytic in $\D$. Then 
$f\in A^p_\alpha$ if and only if there exists a continuous function
$g\in L^p(\D,dA_\alpha)$ such that
\begin{equation}
|f(z)-f(w)|\le\rho(z,w)(g(z)+g(w))
\label{eq2}
\end{equation}
for all $z$ and $w$ in $\D$.
\label{6}
\end{thm}

\begin{proof}
First assume that condition (\ref{eq2}) holds. Then
$$\left|\frac{f(z)-f(w)}{z-w}\right|\le\frac{\rho(z,w)}{|z-w|}(g(z)+g(w))$$
for all $z\not=w$ in $\D$. Fix any $z\in\D$, let $w\to z$, and use
Lemma~\ref{4}. We obtain
$$(1-|z|^2)|f'(z)|\le 2g(z),\qquad z\in\D.$$
Since $g\in L^p(\D,\daa)$, an application of Lemma~\ref{5} shows that
$f\in A^p_\alpha$.

Next assume that $f\in A^p_\alpha$. We are going to produce a continuous
function $g$ that satisfies condition (\ref{eq2}). To this end we fix
a radius $r\in(0,1)$ and consider any two points $z$ and $w$ in $\D$
with $\rho(z,w)<r$. By the fundamental theorem of calculus,
$$f(z)-f(w)=(z-w)\int_0^1f'\bigl(tz+(1-t)w\bigr)\,dt.$$
Since $z$ and $w$ are points in the convex set $D(z,r)$, we see that
$$tz+(1-t)w\in D(z,r)$$
for all $t\in[0,1]$. It follows that
$$|f(z)-f(w)|\le|z-w|\sup\{|f'(u)|:u\in D(z,r)\}.$$
By Lemma~\ref{2}, there exists a positive constant $C$ that only depends
on $r$ such that
$$|f(z)-f(w)|\le\rho(z,w)h(z),$$
where
$$h(z)=C\sup\{(1-|u|^2)|f'(u)|:u\in D(z,r)\},\qquad z\in\D.$$
Obviously, we also have
$$|f(z)-f(w)|\le\rho(z,w)(h(z)+h(w)),$$
provided that $\rho(z,w)<r$.

If $\rho(z,w)\ge r$, then clearly
$$|f(z)-f(w)|\le\frac{\rho(z,w)}r(|f(z)|+|f(w)|).$$
If we define
$$g(z)=\frac{|f(z)|}r+h(z),\qquad z\in\D,$$
then we have
$$|f(z)-f(w)|\le\rho(z,w)(g(z)+g(w))$$
for all $z$ and $w$ in $\D$.

It is clear that the function $g$ above is continuous on $\D$. It remains
for us to show that $g\in L^p(\D,\daa)$. Since $f$ is already in $L^p(\D,\daa)$,
it suffices for us to show that $h\in L^p(\D,\daa)$.

Recall that if $r$ and $R$ are related as in (\ref{eq1}), then $D(z,r)=E(z,R)$.
So we can choose $r'\in(0,1)$ so that $D(z,r')=E(z,2R)$ for all $z\in\D$. By
the triangle inequality for the Bergman metric $\beta$ we have
$E(u,R)\subset E(z,2R)$ whenever $u\in E(z,R)$. Equivalently, we have
$D(u,r)\subset D(z,r')$ whenever $u\in D(z,r)$. It follows from this and
Lemma~\ref{3} that there exists a constant $C>0$ that only depends on $r$
such that
$$h(z)^p\le\frac C{(1-|z|^2)^{2-p}}\int_{D(z,r')}|f'(w)|^p\,dA(w)$$
for all $z\in\D$. If we use $\chi_z$ to denote the characteristic function of
the set $D(z,r')$, then clearly $\chi_z(w)=\chi_w(z)$, and
$$h(z)^p\le\frac C{(1-|z|^2)^{2-p}}\ind|f'(w)|^p\chi_z(w)\,dA(w)$$
for all $z\in\D$. Writing $C_1=(\alpha+1)C/\pi$ and using Fubini's theorem, 
we obtain
\begin{eqnarray*}
\ind h(z)^p\,\daa(z)&\le& C_1\ind(1-|z|^2)^{p+\alpha-2}\,dA(z)\!\!\ind|f'(w)|^p
\chi_z(w)\,dA(w)\\
&=&C_1\ind|f'(w)|^p\,dA(w)\!\!\ind(1-|z|^2)^{p+\alpha-2}\chi_w(z)\,dA(z)\\
&=&C_1\ind|f'(w)|^p\,dA(w)\!\!\int_{D(w,r')}(1-|z|^2)^{p+\alpha-2}\,dA(z).
\end{eqnarray*}
Combining this with Lemma~\ref{2} and the fact that the area of $D(w,r')$ is
comparable to $(1-|w|^2)^2$, we obtain another positive constant $C_2$, which
only depends on $\alpha$ and $r$, such that
$$\ind h(z)^p\,\daa(z)\le C_2\ind|f'(w)|^p(1-|w|^2)^{p+\alpha}\,dA(w).$$
In view of Lemma~\ref{5}, this shows that $h\in L^p(\D,\daa)$ and completes
the proof of Theorem~\ref{6}.
\end{proof}

\begin{thm}
Suppose $p>0$, $\alpha>-1$, and $f$ is analytic in $\D$. Then $f\in A^p_\alpha$
if and only if there exists a continuous function $g\in L^p(\D,\daa)$ such that
\begin{equation}
|f(z)-f(w)|\le\beta(z,w)(g(z)+g(w))
\label{eq3}
\end{equation}
for all $z$ and $w$ in $\D$.
\label{7}
\end{thm}

\begin{proof}
If condition (\ref{eq3}) is satisfied, we divide both sides of (\ref{eq3}) by
$|z-w|$ and, with the help of Lemma~\ref{4}, take the limit as $w\to z$. The
result is
$$(1-|z|^2)|f'(z)|\le 2g(z),\qquad z\in\D.$$
This along with Lemma~\ref{5} shows that $f\in A^p_\alpha$.

If $f\in A^p_\alpha$, then by Theorem~\ref{6}, there exists a continuous function 
$g\in L^p(\D,\daa)$ such that condition (\ref{eq2}) holds. Since 
$\rho(z,w)\le\beta(z,w)$ for all $z$ and $w$ in $\D$, the same function $g$ also 
satisfies condition (\ref{eq3}). This completes the proof of Theorem~\ref{7}.
\end{proof}

\begin{thm}
Suppose $p>0$, $\alpha>-1$, and $f$ is analytic in $\D$. Then $f\in A^p_\alpha$
if and only if there exists a continuous function $g\in L^p(\D,dA_{p+\alpha})$ 
such that
\begin{equation}
|f(z)-f(w)|\le|z-w|(g(z)+g(w))
\label{eq4}
\end{equation}
for all $z$ and $w$ in $\D$.
\label{8}
\end{thm}

\begin{proof}
If condition (\ref{eq4}) holds, we divide both sides of (\ref{eq4}) by $|z-w|$
and take the limit as $w\to z$. The result is $|f'(z)|\le 2g(z)$, so that
$$(1-|z|^2)|f'(z)|\le 2(1-|z|^2)g(z)$$
for all $z\in\D$. Since $g\in L^p(\D,dA_{p+\alpha})$, we see that the
function 
$$(1-|z|^2)f'(z)$$
belongs to $L^p(\D,\daa)$, which, according to Lemma~\ref{5}, means that 
$f\in A^p_\alpha$.

If $f\in A^p_\alpha$, then by Theorem~\ref{6}, there exists a continuous
function $h$ in $L^p(\D,\daa)$ such that
$$|f(z)-f(w)|\le\rho(z,w)(h(z)+h(w))$$
for all $z$ and $w$ in $\D$. Rewrite this as
$$|f(z)-f(w)|\le|z-w|\left[\frac{h(z)}{|1-\overline zw|}+
\frac{h(w)}{|1-\overline zw|}\right],$$
and apply the triangle inequalities
$$|1-\overline zw|\ge1-|z|,\quad |1-\overline zw|\ge1-|w|.$$
We obtain
$$|f(z)-f(w)|\le|z-w|(g(z)+g(w)),\qquad z,w\in\D,$$
where
$$g(z)=\frac{h(z)}{1-|z|}\le\frac{2h(z)}{1-|z|^2}.$$
Since $h\in L^p(\D,\daa)$, we have $g\in L^p(\D,dA_{p+\alpha})$. This
completes the proof of Theorem~\ref{8}.
\end{proof}

\section{Lifting functions from the disk to the bidisk}

Let $\D^2=\D\times\D$ denote the bidisk in $\C^2$ and let $H(\D^2)$ denote
the space of all holomorphic functions in $\D^2$. Similarly, $H(\D)$ is
the space of all analytic functions in $\D$. For $p>0$ and
$\alpha>-1$ we define $A^p_\alpha(\D^2)$ as the space of functions
$f\in H(\D^2)$ such that
$$\ind\ind|f(z,w)|^p\,dA_\alpha(z)\,dA_\alpha(w)<\infty.$$
These are also called weighted Bergman spaces.

In this section we present an application of our main theorem to the
problem of lifting analytic functions from the unit disk to the bidisk.
Thus we consider the symmetric lifting operator
$$L: H(\D)\to H(\D^2)$$
defined by
$$L(f)(z,w)=\frac{f(z)-f(w)}{z-w}.$$
We will also need the associated diagonal operator
$$\Delta: H(\D^2)\to H(\D)$$
which is defined by
$$\Delta(f)(z)=f(z,z).$$
The action of the diagonal operator on Hardy and Bergman spaces of the
polydisk has been studied by several authors. See \cite{DS75}
\cite{HO}\cite{FAS}\cite{SZ}. In particular, the diagonal operator has 
the following property.

\begin{lemma}
Suppose $p>0$ and $\alpha>-1$. Then the operator $\Delta$ maps 
$A^p_\alpha(\D^2)$ boundedly onto $A^p_{2(\alpha+1)}(\D)$.
\label{9}
\end{lemma}

\begin{proof}
See \cite{FAS} or \cite{SZ}.
\end{proof}

The following standard estimate will be needed in the proof of our lifting
theorems.

\begin{lemma}
Suppose $s>-1$, $t$ is real, and
$$I(z)=\ind\frac{(1-|w|^2)^s\,dA(z)}{|1-\overline zw|^{2+s+t}},
\qquad z\in\D.$$
Then $I(z)$ is bounded in $\D$ whenever $t<0$; and $I(z)$ is bounded by
$(1-|z|^2)^{-t}$ whenever $t>0$.
\label{10}
\end{lemma}

\begin{proof}
See \cite{HKZ} or \cite{Z}.
\end{proof}

We now obtain the first lifting theorem.

\begin{thm}
Suppose $\alpha>-1$ and $0<p<\alpha+2$. Then the symmetric lifting 
operator $L$ maps $A^p_\alpha(\D)$  boundedly into
$A^p_\alpha(\D^2)$. Moreover, this is no longer true when $p>\alpha+2$.
\label{11}
\end{thm}

\begin{proof}
Given $f\in A^p_\alpha(\D)$, we apply Theorem~\ref{6} to find a function
$g\in L^p(\D,\daa)$ such that condition (\ref{eq2}) holds. Then there 
exists a constant $C=C_p$ such that
\begin{equation}
\left|L(f)(z,w)\right|^p\le C\left[\frac{g(z)^p}{|1-\overline zw|^p}
+\frac{g(w)^p}{|1-\overline zw|^p}\right].
\label{eq5}
\end{equation}
It follows that
$$\ind\ind\left|L(f)(z,w)\right|^p\,dA_\alpha(z)\,dA_\alpha(w)\le 2C
\ind g(z)^p\,dA_\alpha(z)\ind\frac{dA_\alpha(w)}{|1-\overline zw|^p}.$$
When $p<2+\alpha$, an application of Lemma~\ref{10} shows that there
exists another constant $C>0$ such that
$$\ind\ind\left|L(f)(z,w)\right|^p\,dA_\alpha(z)\,dA_\alpha(w)\le C
\ind g(z)^p\,dA_\alpha(z).$$
This shows that $L$ maps $A^p_\alpha(\D)$ into $A^p_\alpha(\D^2)$. A
standard argument based on the closed graph theorem then shows that
the operator
$$L:A^p_\alpha(\D)\to A^p_\alpha(\D^2)$$
must be bounded.

On the other hand, let us suppose that
$$L:A^p_\alpha(\D)\to A^p_\alpha(\D^2)$$
is bounded. Then by Lemma~\ref{9}, the operator $D=\Delta\circ L$ maps
$A^p_\alpha(\D)$ boundedly into $A^p_{2(\alpha+1)}(\D)$. It is easy to
see that $Df=f'$. Therefore, $f\in A^p_\alpha(\D)$ would imply that
$$\ind\left|(1-|z|^2)f'(z)\right|^p\,dA_{2(\alpha+1)-p}(z)<\infty,$$
which, according to Lemma~\ref{5}, is a condition that is strictly 
stronger than $f\in A^p_\alpha(\D)$ when $p>2+\alpha$. So our lifting
theorem cannot possibly hold for $p>2+\alpha$.
\end{proof}

We mention that the case $p=2+\alpha$ is not covered by the above
result. When $\alpha=0$, we can show by Taylor expansion that the
operator $L$ does not map $A^2(\D)$ into $A^2(\D^2)$ (these are the
unweighted Bergman spaces). In fact, if 
$$f(z)=\sum_{k=0}^\infty a_kz^k$$
is a function in $A^2(\D)$, then
$$\ind|f(z)|^2\,dA(z)=\sum_{k=0}^\infty\frac{|a_k|^2}{k+1}.$$
On the other hand,
$$L(f)(z,w)=\sum_{k=0}^\infty a_k\frac{z^k-w^k}{z-w}
=\sum_{k=1}^\infty a_k\sum_{i+j=k-1}z^iw^j,$$
and for $k\not=m$, the homogeneous polynomials
$$\sum_{i+j=k-1}z^iw^j,\quad \sum_{i+j=m-1}z^iw^j,$$
are orthogonal with respect to the measure $dA(z)dA(w)$ on $\D^2$. 
So the integral
$$I=\ind\ind\left|L(f)(z,w)\right|^2\,dA(z)\,dA(w)$$
can be computed as follows.
\begin{eqnarray*}
I&=&\sum_{k=1}^\infty|a_k|^2\ind\ind\left|\sum_{i+j=k-1}
z^iw^j\right|^2\,dA(z)\,dA(w)\\
&=&\sum_{k=1}^\infty|a_k|^2\sum_{i+j=k-1}\frac1{(i+1)(j+1)}\\
&=&\sum_{k=1}^\infty\frac{|a_k|^2}{k+1}\sum_{i+j=k-1}\left[
\frac1{i+1}+\frac1{j+1}\right]\\
&=&2\sum_{k=1}^\infty\frac{|a_k|^2}{k+1}\sum_{j=0}^{k-1}
\frac1{j+1}\\
&\sim&\sum_{k=1}^\infty\frac{|a_k|^2}{k+1}\log(k+1)\\
&\sim&\ind|f(z)|^2\log\frac1{1-|z|^2}\,dA(z).
\end{eqnarray*}
This shows that the integral $I$ is not necessarily finite, so the 
symmetric lifting operator $L$ does not map $A^2(\D)$ into $A^2(\D^2)$.

The following result tells us what happens when $p>\alpha+2$.

\begin{thm}
Suppose $\alpha>-1$, $p>\alpha+2$, and $\beta$ is determined by
$$2(\beta+1)=p+\alpha.$$
Then the symmetric lifting operator $L$ maps
$A^p_\alpha(\D)$ boundedly into $A^p_\beta(\D^2)$.
\label{12}
\end{thm}

\begin{proof}
Given $f\in A^p_\alpha(\D)$, we once again appeal to Theorem~\ref{6}
to obtain a function $g\in L^p(\D,dA_\alpha)$ such that 
condition (\ref{eq5}) holds. We then have
$$\ind\ind\left|L(f)(z,w)\right|^p\,dA_\beta(z)\,dA_\beta(w)\le 2C
\ind g(z)^p\,dA_\beta(z)\ind\frac{dA_\beta(w)}{|1-\overline zw|^p}.$$
Since $p+\alpha=2(\beta+1)$ and $p>\alpha+2$, we must have $\beta>-1$
and $\beta>\alpha$. We write the inner integral above as
$$(\beta+1)\ind\frac{(1-|w|^2)^\beta\,dA(w)}
{|1-\overline zw|^{2+\beta+(\beta-\alpha)}}$$
and apply Lemma~\ref{10} to find another positive constant $C$ such that
$$\ind\ind\left|L(f)(z,w)\right|^p\,dA_\beta(z)\,dA_\beta(w)\le C
\ind g(z)^p\,dA_\alpha(z).$$
This along with the closed graph theorem proves the desired result.
\end{proof}

Once again, with the help of Lemmas~\ref{9} and \ref{5}, we can check
that the lifting effect of $L$ guaranteed by Theorem~\ref{12} is best 
possible.

\section{Generalization to the Unit Ball}

In this section we explain how to generalize our main result to the
context of the unit ball in $\cn$. Thus we let 
$$\bn=\{z\in\cn:|z|<1\}$$ 
denote the open unit ball in $\cn$. For $\alpha>-1$ let
$$dv_\alpha(z)=c_\alpha(1-|z|^2)^\alpha\,dv(z),$$
where $dv$ is normalized volume measure on $\bn$ and $c_\alpha$ is a
positive normalizing constant so that $v_\alpha(\bn)=1$.

For $p>0$ and $\alpha>-1$ let
$$A^p_\alpha(\bn)=L^p(\bn,dv_\alpha)\cap H(\bn)$$
denote the weighted Bergman spaces, where $H(\bn)$ is the space of all 
holomorphic functions in $\bn$. See \cite{Z} for basic properties of 
these spaces.

For any $a\in\bn$ there exists a biholomorphic map $\varphi_a$
on $\bn$ such that $\varphi_a(0)=a$ and $\varphi^{-1}_a=\varphi_a$.
These are sometimes called symmetrices (or involutive automorphisms)
of $\bn$. Explicit formulas are available for $\varphi_a$; see \cite{R}
or \cite{Z}.

It is well known that the Bergman metric on $\bn$ induces the following
distance:
$$\beta(z,w)=\frac12\log\frac{1+|\varphi_z(w)|}{1-|\varphi_z(w)|}.$$
It follows that
$$\rho(z,w)=|\varphi_z(w)|$$
is also a distance function on $\bn$. We shall also call $\rho$ the
pseudo-hyperbolic metric on $\bn$. The Euclidean metric on $\bn$ is
still denoted by $|z-w|$.

\begin{thm}
Suppose $p>0$, $\alpha>-1$, and $f\in H(\bn)$. Then the following
conditions are equivalent.
\begin{itemize}
\item[(a)] $f$ belongs to $A^p_\alpha(\bn)$.
\item[(b)] There exists a continuous function $g\in L^p(\bn,dv_\alpha)$
such that
$$|f(z)-f(w)|\le\rho(z,w)(g(z)+g(w))$$
for all $z$ and $w$ in $\bn$.
\item[(c)] There exists a continuous function $g\in L^p(\bn,dv_\alpha)$
such that
$$|f(z)-f(w)|\le\beta(z,w)(g(z)+g(w))$$
for all $z$ and $w$ in $\bn$.
\item[(d)] There exists a continuous function $g\in L^p(\bn,dv_{p+\alpha})$
such that
$$|f(z)-f(w)|\le|z-w|(g(z)+g(w))$$
for all $z$ and $w$ in $\bn$.
\end{itemize}
\label{13}
\end{thm}

The proof follows the same ideas used in the proof of Theorems~\ref{6},
\ref{7}, and \ref{8}. Lemmas~\ref{2} holds for the ball without any
change; see \cite{Z}. 

The only change needed in Lemma~\ref{3} is the expenent of $1-|z|^2$. In 
the context of $\bn$, it should be $(1-|z|^2)^{n+1}$ instead of 
$(1-|z|^2)^2$; see \cite{Z}. Similarly, the volume of $D(z,r)$ (or $E(z,R)$) 
is comparable to $(1-|z|^2)^{n+1}$ whenever $r$ (or $R$) is fixed.

Lemmas~\ref{4} and \ref{5} need to be modified substantially before they 
can be used. These are the contents of the next two lemmas. But first,
we recall three useful notations of differentiation in $\bn$.

Given $f\in H(\bn)$, we write
$$Rf(z)=\sum_{k=1}^nz_k\frac{\partial f}{\partial z_k}(z),$$
and call it the radial derivative of $f$ at $z$. In fact, $Rf(z)$ is
the directional derivative of $f$ at $z$ in the radial direction (that
is, the direction in $z$):
$$Rf(z)=\lim_{t\to1}\frac{f(tz)-f(z)}t,$$
where $t$ is a scalar.

The complex gradient of $f$ at $z$ is defined by
$$|\nabla f(z)|=\left[\sum_{k=1}^n\left|\frac{\partial f}{\partial z_k}(z)
\right|^2\right]^{1/2}.$$
And the invariant complex gradient of $f$ at $z$ is given by
$$|\widetilde\nabla f(z)|=|\nabla(f\circ\varphi_z)(0)|,\qquad z\in\bn.$$

We can now state the analogs of Lemmas~\ref{4} and \ref{5} that will suit
our needs.

\begin{lemma}
Suppose $z\in\bn$ and $w=tz$, where $t$ is a scalar. Then
$$\lim_{w\to z}\frac{\rho(w,z)}{|z-w|}=\lim_{w\to z}\frac{\beta(z,w)}
{|z-w|}=\frac1{1-|z|^2}.$$
\label{14}
\end{lemma}

\begin{proof}
This follows from the explicit formulas for $\varphi_a$ given in \cite{R}
and \cite{Z}.
\end{proof}

\begin{lemma}
Suppose $p>0$, $\alpha>-1$, and $f\in H(\bn)$. Then the following
conditions are equivalent.
\begin{itemize}
\item[(a)] The function $f$ is in $A^p_\alpha(\bn)$.
\item[(b)] The function $(1-|z|^2)Rf(z)$ is in $L^p(\bn,dv_\alpha)$.
\item[(c)] The function $(1-|z|^2)|\nabla f(z)|$ is in $L^p(\bn,dv_\alpha)$.
\item[(d)] The function $|\widetilde\nabla f(z)|$ is in $L^p(\bn,dv_\alpha)$.
\end{itemize}
\label{15}
\end{lemma}

\begin{proof}
See \cite{Z}.
\end{proof}

We can now outline the proof of Theorem~\ref{13}.

First assume that there exists a continuous function 
$g\in L^p(\bn,dv_\alpha)$ such that
\begin{equation}
|f(z)-f(w)|\le\rho(z,w)(g(z)+g(w))
\label{eq6}
\end{equation}
for all $z$ and $w$ in $\bn$. We then fix $z$ in $\bn$ and let $w=tz$,
where $t$ is a scalar. Then
$$\frac{|f(z)-f(w)|}{|z-w|}\le\frac{\rho(z,w)}{|z-w|}(g(z)+g(w))$$
for all $z\not=w$ in $\bn$. Let $w$ approach $z$ in the radial direction
and apply Lemma~\ref{14}. We obtain
$$(1-|z|^2)|Rf(z)|\le 2g(z)$$
for all $z\in\bn$. According to Lemma~\ref{15}, this is the same as
$f\in A^p_\alpha(\bn)$.

On the other hand, for any holomorphic function $f$ in $\bn$ and any
$z\in\bn$, we have
$$f(z)-f(0)=\int_0^1\left[\sum_{k=1}^nz_k\frac{\partial f}{\partial z_k}
(tz)\right]\,dt.$$
It follows that for $\rho(z,0)<r$, where $r\in(0,1)$ is any fixed
radius in the pseudo-hyperbolic metric, we have
$$|f(z)-f(0)|\le|z|\,|\sup\{|\nabla f(w)|:w\in D(0,r)\}.$$
It is easy to see that in the relatively compact set $D(0,r)$ the
Euclidean metric is comparable to the pseudo-hyperbolic metric (as well
as the Bergman metric $\beta$). It is also easy to see that $|\nabla f(w)|$
is comparable to $|\widetilde\nabla f(w)|$ in the relatively compact
set $D(0,r)$. So we can find a constant $C>0$, that depends on $r$ but
not on $f$, such that
$$|f(z)-f(0)|\le C\rho(z,0)\sup\{|\widetilde\nabla f(w)|:w\in D(0,r)\}$$
for all $z\in D(0,r)$. Replace $f$ by $f\circ\varphi_w$, then replace $z$
by $\varphi_w(z)$, and use the M\"obius invariance of the pseudo-hyperbolic 
metric and the invariant gradient. We obtain
$$|f(z)-f(w)|\le C\rho(z,w)\sup\{|\widetilde\nabla f(u)|:u\in D(z,r)\}$$
for all $z$ and $w$ in $\bn$ with $\rho(z,w)<r$. Let
$$g(z)=\frac{|f(z)|}r+C\sup\{|\widetilde\nabla f(u)|:u\in D(z,r)\}.$$
Then it is clear that condition (\ref{eq6}) is satisfied, and just like
in the proof of Theorem~\ref{6}, $g\in L^p(\bn,dv_\alpha)$.

So conditions (a) and (b) are equivalent in Theorem~\ref{13}.

If condition (b) holds in Theorem~\ref{13}, then condition (c) holds for
the same function $g$, because $\rho\le\beta$. If condition (c) holds,
then an application of Lemma~\ref{14} shows that
$$(1-|z|^2)|Rf(z)|\le 2g(z)$$
for all $z\in\bn$, which, according to Lemma~\ref{15}, implies that 
condition (a). Therefore, conditions (a), (b), and (c) are all equivalent.

Now let us assume that condition (d) holds, so there exists a continuous
function $g$ in $L^p(\bn,dv_{p+\alpha})$ such that
\begin{equation}
|f(z)-f(w)|\le|z-w|(g(z)+g(w))
\label{eq7}
\end{equation}
for all $z$ and $w$ in $\bn$. Rewrite this as
$$\frac{|f(z)-f(w)|}{|z-w|}\le g(z)+g(w)$$
and let $w$ approach $z$ in the direction of a complex coordinate axis.
We obtain
$$\left|\frac{\partial f}{\partial z_k}(z)\right|\le 2g(z),
\qquad 1\le k\le n,$$
so that
$$|\nabla f(z)|\le2\sqrt n\,g(z)$$
for all $z\in\bn$. This together with the assumption that $g\in L^p(\bn,
dv_{p+\alpha})$ shows that the function $(1-|z|^2)|\nabla f(z)|$ is
in $L^p(\bn,dv_\alpha)$. In view of Lemma~\ref{15}, this is the same
as $f\in A^p_\alpha(\bn)$. So condition (d) implies (a) in Theorem~\ref{13}.

Finally let us assume that condition (b) holds in Theorem~\ref{13}. It 
follows from the well-known identity (see \cite{R} or \cite{Z})
$$1-|\varphi_z(w)|^2=\frac{(1-|z|^2)(1-|w|^2)}{|1-\langle z,w\rangle|^2}$$
that
$$\frac{\rho(z,w)^2}{|z-w|^2}=\frac{|z-w|^2+|\langle z,w\rangle|^2-
|z|^2|w|^2}{|z-w|^2|1-\langle z,w\rangle|^2}.$$
By the triangle inequality for the natural inner product in $\cn$, we
always have
$$|\langle z,w\rangle|^2\le|z|^2|w|^2.$$
We deduce that
$$\rho(z,w)\le\frac{|z-w|}{|1-\langle z,w\rangle|}$$
for all $z$ and $w$ in $\bn$. Just like in the proof of Theorem~\ref{8},
this along with condition (b) implies the existence of a continuous function 
$g$ in the space $L^p(\bn,dv_{p+\alpha})$ such that condition (\ref{eq7}) holds.
So condition (b) implies (d), and the proof of Theorem~\ref{13} is complete.

Although we have not checked whether the function
$$d(z,w)=\frac{|z-w|}{|1-\langle z,w\rangle|},\qquad z,w\in\bn,$$
defines a metric on $\bn$, it is clear by now that the pseudo-hyperbolic
metric $\rho$ used in Theorem~\ref{13} can be replaced by $d$. In other
words, a holomorphic function $f$ in $\bn$ belongs to the Bergman space
$A^p_\alpha(\bn)$ if and only if there exists a continuous function
$g\in L^p(\bn,dv_\alpha)$ such that
$$|f(z)-f(w)|\le d(z,w)(g(z)+g(w))$$
for all $z$ and $w$ in $\bn$.

\end{document}